\documentclass{ifacconf}

\usepackage{graphicx}      
\usepackage{natbib}        
\usepackage{url} 
\usepackage{amsmath,amssymb,amsfonts}
\usepackage{algorithmic}
\usepackage{textcomp}
\usepackage{xcolor}
\usepackage[short]{optidef}
\usepackage{tikz}
\usetikzlibrary{arrows,arrows.meta,shapes,positioning,shadows,trees}
\usepackage{siunitx}
\usepackage{lipsum}
\usepackage{systeme}
\usepackage{siunitx}

\NewDocumentCommand{\vartwo}{ >{\SplitArgument{2}{,}}m }{ \finalvartwo#1 }
\NewDocumentCommand{\finalvartwo}{mmm}{ \ensuremath{#1_\text{#2}} }

\NewDocumentCommand{\varthree}{ >{\SplitArgument{3}{,}}m }{ \finalvarthree#1 }
\NewDocumentCommand{\finalvarthree}{mmmm}{ \ensuremath{#1_{\text{#2},#3}} }
\begin{document}

\begin{frontmatter}
{\color{blue}
This work has been submitted to IFAC for possible publication 

Submitted to 10th Vienna International Conference
on Mathematical Modelling}
\end{frontmatter}

\begin{frontmatter}

\title{Mathematical Modeling for Holistic\\Convex Optimization of Hybrid Trains} 

\author{Rabee Jibrin, Stuart Hillmansen, Clive Roberts}

\address{Department of Electronic, Electrical and Systems Engineering, University of Birmingham, Birmingham, United Kingdom\\ \textup{e-mail: rxj956 , s.hillmansen , c.roberts.20 @bham.ac.uk}}

\begin{abstract}                
We look into modeling fuel cell hybrid trains for the purpose of optimizing their operation using convex optimization. Models and constraints necessary to form a physically feasible yet convex problem are reviewed. This effort is described as holistic due to the broad consideration of train speed, energy management system, and battery thermals. The minimized objective is hydrogen fuel consumption for a given target journey time. A novel battery thermal model is proposed to aid with battery thermal management and thus preserve battery lifetime. All models are derived in the space-domain which along constraint relaxations guarantee a convex optimization problem. First-principle knowledge and real-world data justify the suitableness of the proposed models for the intended optimization problem.
\end{abstract}


\end{frontmatter}

\section{Introduction}

\subsection{Motivation}

Hydrogen fuel cell hybrid trains are expected to play a key role in decarbonizing the railways owing to their lack of harmful emissions at point-of-use and adequate driving range; however, their total cost of ownership is projected to be higher than incumbent diesel trains primarily due to the higher cost of hydrogen fuel in comparison to diesel fuel \citep{RN200}. We aim at reducing hydrogen fuel consumption by optimizing train operation. Models for a convex optimization problem are sought after in order to alleviate computational concerns.

\subsection{Background}

Train speed optimization has been researched extensively owing to the large contribution of traction power towards rail energy consumption \citep{RN707}. More recently, \cite{RN876} brought attention to the energy management system (EMS) of fuel cell hybrid trains by holding the 2019 IEEE VTS Motor Vehicles Challenge. The EMS determines power distribution among a hybrid vehicle's power-sources and is thus a vital determinant of efficiency. \cite{RN318} reviewed an extensive body of literature for fuel cell hybrid EMS. Simulation results by \cite{RN565} suggest that optimization-based algorithms outperform their rule-based counterparts which motivates our current focus on the former.

While the aforementioned address either speed or EMS separately, some works have attempted to optimize both within a single optimization problem (concurrently) in order to achieve better solution optimality by embedding knowledge of the dynamic coupling between both trajectories, e.g., dynamic programming \citep{RN857}, indirect optimal control \citep{RN613}, integer programming \citep{RN924}, and relaxed convex optimisation \citep{RN886,mylatest}.

The high capital cost of traction batteries has also motivated many to consider penalizing \citep{RN380} or bounding \citep{RN384} battery degradation, though strictly within the EMS problem setting---speed is optimized beforehand separately. The semi-empirical battery degradation model presented by \cite{RN370} as a function of temperature, state-of-charge, and C-rate, is the most often used. A common assumption among optimization-based algorithms that consider battery degradation is an active cooling system that maintains a constant battery temperature which simplifies the degradation model to static temperature. This simplification can lead to unexpected battery degradation when subject to non-ideal thermal management in the real-world \citep{Filippi}. Therefore, dropping the static temperature assumption could further benefit battery lifetime, especially in light of experimental results that designate elevated temperatures as the leading cause of battery degradation \citep{RN386}. Moreover, including thermal constraints during speed planning can limit the reliance on the active cooling system and thus reduce its parasitic energy draw \citep{RN926}. Algorithms that did consider battery temperature as a bounded dynamic state have done so strictly within the EMS problem setting and often at great computational cost, e.g., genetic algorithm \citep{RN296}, dynamic programming \citep{RN669}, and relaxed convex optimization \citep{RN808}.

\subsection{Contribution and Outline}

Literature lacks a method to concurrently optimize hybrid vehicle speed and EMS while considering battery thermal constraints. The high predictability of railway environments promises substantial returns for such an elaborate and holistic optimization approach. This paper gathers the models necessary to form a convex optimization problem for the aforementioned goal. Furthermore, a novel thermal model for the battery is proposed. Future publications will showcase these models within a realistic optimization case study, though preliminary results by the authors can be found in \citep{mylatest}.

Section 2 introduces the train's longitudinal dynamics, section 3 covers the powertrain's models, and section 4 uses these models to formulate the optimization problem.

\section{Longitudinal Dynamics}

\subsection{Choice of Modeling Domain}
Common among model-based optimization for dynamic systems is to model the system in the time-domain, i.e., the model predicts system state after a temporal interval of $\vartwo{\Delta,t}$ seconds. However, a complication from optimizing vehicle speed in the time-domain is interpolating track information, e.g., gradient, when the physical location for a given temporal interval is dependent on the optimized speed and thus unknown \textit{a priori}. This can be addressed by relying on historical speed data to predict location against time however significant location errors could accumulate over a long journey. Alternatively, more sophisticated methods such as the pseudospectral method can be used at a great computational cost \citep{RN707}.

Instead of the often used time-domain, the current problem setting lends itself more readily to the space-domain, i.e., the model predicts system state after a spatial interval of $\vartwo{\Delta,s}$ meters longitudinally along the track. As such, one can accurately retrieve track information for any interval by directly referring to its respective location in space. Herein, we formulate the models in the discrete space-domain with zero-order hold between intervals. The entire journey's longitudinal space is divided into a grid of $N$ intervals.

\subsection{Train Longitudinal Speed}

\cite{RN728} assume the train as a point mass $m$ with an equivalent inertial mass $\vartwo{m,eq}$. The train's longitudinal speed $v$ is influenced by traction motor force $\vartwo{F,m}$, mechanical brakes force $\vartwo{F,brk}$, and the external forces acting on the train $\vartwo{F,ext}$ which is the summation of the Davis Equation $a+bv_i+cv_i^2$ and gravitational pull $mg\sin(\theta_i)$. To predict speed after a single spatial interval, construct
\begin{equation}\label{eq:kinetic_energy_1}
\frac{1}{2}\vartwo{m,eq}v_{i+1}^2 = \frac{1}{2}\vartwo{m,eq}v_{i}^2 + (\varthree{F,m,i} + \varthree{F,brk,i})\Delta_{\text{s},i} - \varthree{F,ext,i}\Delta_{\text{s},i}
\end{equation}

using the definition of kinetic energy $\vartwo{E,k.e.}=1/2 \vartwo{m,eq} v^2$, the definition of mechanical work $\vartwo{E,work}=F\Delta_\text{s}$, and the principle of energy conservation. Equation \eqref{eq:kinetic_energy_1} is nonlinear in $v$ but can be linearized by substituting the quadratic terms $v^2$ with $z$ and keeping the non-quadratic terms $v$ unchanged, namely
\begin{equation}\label{eq:kinetic_energy_2}
\frac{1}{2}\vartwo{m,eq}z_{i+1} = \frac{1}{2}\vartwo{m,eq}z_{i} + (\varthree{F,m,i} + \varthree{F,brk,i})\Delta_{\text{s},i} - \varthree{F,ext,i}\Delta_{\text{s},i}
\end{equation}
and
\begin{equation}\label{eq:f_ext}
\varthree{F,ext,i} = a+bv_i+cz_i + mg\sin(\theta_i). 
\end{equation}

The linear model \eqref{eq:kinetic_energy_2} relies on both $v$ and $z$ to define train speed and thus requires the non-convex equality constraint $v^2=z$ to hold true which is subsequently relaxed into the convex inequality
\begin{equation}\label{eq:relaxed_v_2}
v^2 \leq z.
\end{equation}

\subsection{Journey Target Time}

Total journey time is expressed as summation of time required for all intervals $\sum_{i=1}^N \vartwo{\Delta,s}/v_i$ but is non-linear in $v$. This expression can be replaced by the linear 
\begin{equation}
\sum_{i=1}^N \vartwo{\Delta,s}\lambda_{v,i}
\end{equation} 
when used along the auxiliary non-convex equality $\lambda_v = 1/v$ which is then relaxed into the convex inequality
\begin{equation}\label{eq:lambda_v}
\lambda_{v} \geq 1/v
\end{equation}
for $v,\lambda_v > 0$ \citep{boyd2004convex}. Section 4 explains how the strict positivity constraints imposed on speed have a negligible impact on solution optimality and how the relaxed inequalities \eqref{eq:relaxed_v_2} and \eqref{eq:lambda_v} hold with equality at the optimal solution.

\section{Powertrain Modeling}

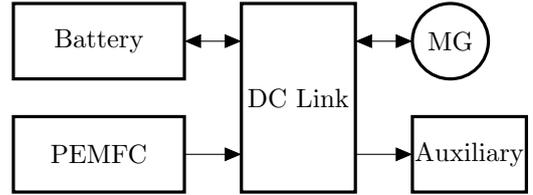
\begin{figure}[!t]
\centering
  \begin{tikzpicture}
  \draw[black, very thick] (-0.25,0) rectangle (2,1) node[pos=0.5] {Battery};
  \draw[black, very thick] (-0.25,-1.5) rectangle (2,-0.5) node[pos=0.5] {PEMFC};
  \draw[black, very thick] (2.75,-1.5) rectangle (4.25,1) node[pos=0.5, align=center] {DC Link};
  \draw[black, very thick] (5.5,0.5) circle (0.5cm) node {MG};
  \draw[black, very thick] (5,-1.5) rectangle (6.5,-0.5) node[pos=0.5] {Auxiliary};
  \draw[>=triangle 45, <->] (2,0.5) -- (2.75,0.5); 
  \draw[>=triangle 45, ->] (2,-1) -- (2.75,-1); 
  \draw[>=triangle 45, <->] (4.25,0.5) -- (5,0.5); 
  \draw[>=triangle 45, ->] (4.25,-1) -- (5,-1); 
  \end{tikzpicture}
\caption{Fuel cell series hybrid architecture. Arrows depict feasible directions of electric power flow.}
\label{fig:fc_powertrain}
\end{figure}

Figure \ref{fig:fc_powertrain} depicts the powertrain considered, a polymer electrolyte membrane fuel cell (PEMFC) in a series hybrid configuration with a lithium-ion battery. The components considered herein are the battery, fuel cell, motor-generator (MG), and vehicle auxiliary loads. The term motor is used interchangeably with motor-generator. The following subsections present the models and constraints for each component. Without loss of generality, repeated components are aggregated and modeled as a single bigger component, optimized as the newly formed single big component, after which the optimized solution is divided equally upon the actual individual instances of that component, e.g., the traction motors are modeled and optimized as a single big motor acting on the point mass.

\subsection{Traction Motor}\label{sec:motor}

\begin{figure}[t!]
\centering
\includegraphics[width=8.4cm,keepaspectratio]{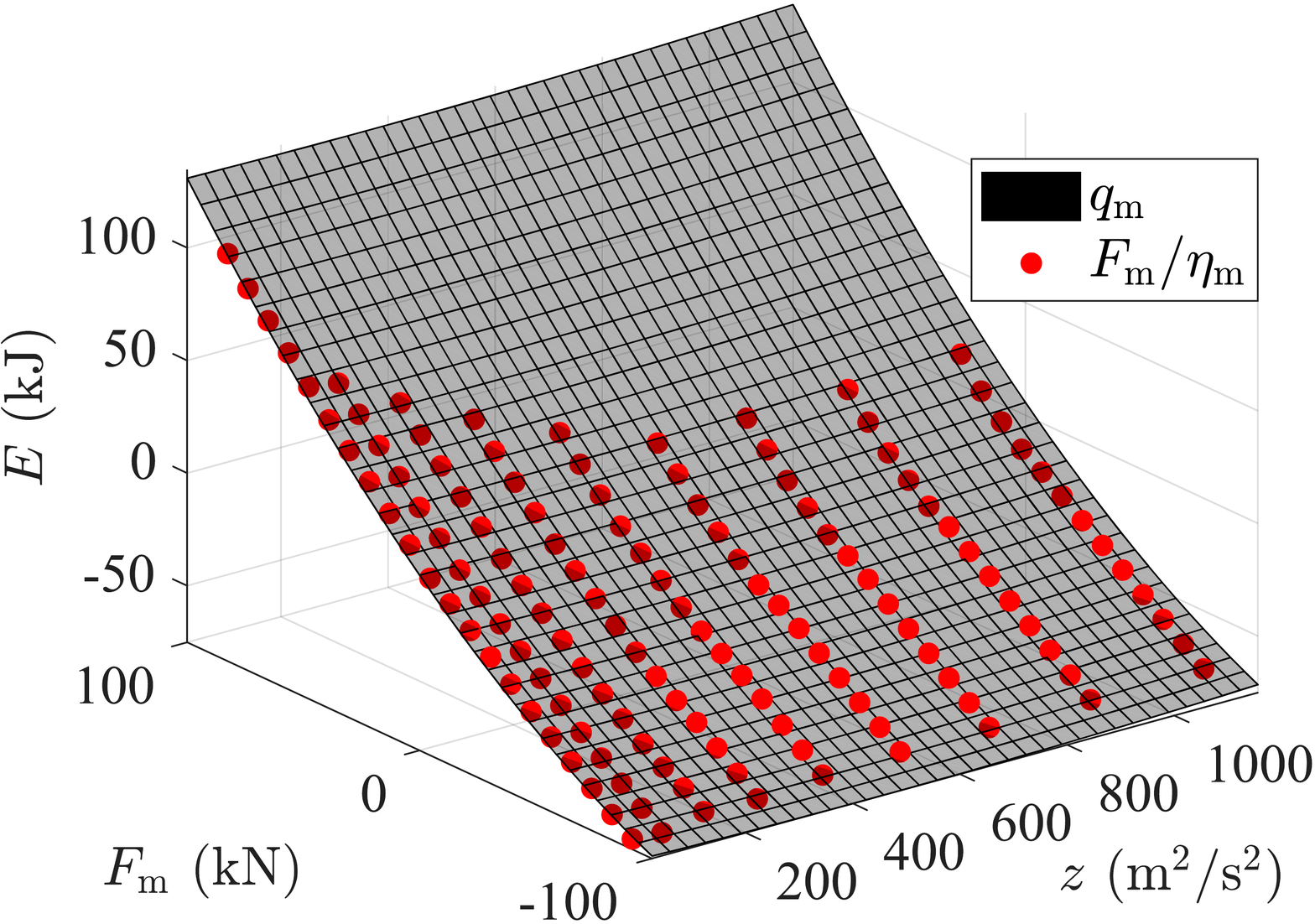}
\caption{Motor data from \cite{RN879}. $\vartwo{\Delta,s}=1$.}
\label{fig:motor_consumption}
\end{figure}

The electric power flow in Fig. \ref{fig:fc_powertrain} is described by 
\begin{equation}\label{eq:power_balance}
\vartwo{P,m}/\vartwo{\eta,m}(\vartwo{P,m}) + \vartwo{P,aux} = \vartwo{P,fc} + \vartwo{P,batt},
\end{equation}
where $\vartwo{P,m}$ is motor mechanical power, $\eta_m(\vartwo{P,m})$ is motor efficiency and thus $\vartwo{P,m}/\eta_m(\vartwo{P,m})$ is electric power at motor terminals, $\vartwo{P,aux}$ is auxiliary load, $\vartwo{P,fc}$ is fuel cell electric power output, and $\vartwo{P,batt}$ is battery electric power output.

The power balance expression \eqref{eq:power_balance} requires the non-convex constraint $\vartwo{P,m}=\vartwo{F,m}v$ to hold true in order to use it in conjunction with the speed model \eqref{eq:kinetic_energy_2}. To resolve this non-convexity, start by dividing \eqref{eq:power_balance} by $v$ to yield
\begin{equation}\label{eq:force_balance}
\vartwo{F,m}/\vartwo{\eta,m}(\vartwo{F,m},z) + \vartwo{P,aux}\lambda_v = \vartwo{F,fc} + \vartwo{F,batt},
\end{equation}
where motor efficiency is defined as $\vartwo{\eta,m}(\vartwo{F,m},z)$ instead of $\vartwo{\eta,m}(\vartwo{P,m})$, recall $P=F\sqrt{z}$. The alternative model \eqref{eq:force_balance} expresses energy flow per longitudinal meter traveled, recall $\vartwo{E,work}=F\vartwo{\Delta,s}$ and $F=P\lambda_v$. The forces $\vartwo{F,fc}$ and $\vartwo{F,batt}$ are fictitious but numerically represent each component's energy contribution. Since \eqref{eq:force_balance} is directly written in terms of $\vartwo{F,m}$ the non-convex equality constraint $\vartwo{P,m}=\vartwo{F,m}v$ is no longer necessary.

The equality \eqref{eq:force_balance} is non-convex due to the non-linearity in $\vartwo{F,m}/\vartwo{\eta,m}(\vartwo{F,m},z)$. Moreover, motor efficiency, $\vartwo{\eta,m}$, is typically a discrete look-up table rather than a smooth function. \cite{RN862} accurately approximated $\vartwo{F,m}/\vartwo{\eta,m}(\vartwo{F,m},z)$ with the convex quadratic polynomial $\vartwo{q,m}(\vartwo{F,m},z):=p_{00}+p_{10}z+p_{01}\vartwo{F,m}+p_{11}\vartwo{F,m}v+p_{20}z^2+p_{02}\vartwo{F,m}^2$, as shown in Fig. \ref{fig:motor_consumption}, which can be used to relax \eqref{eq:force_balance} into the convex inequality
\begin{equation}\label{eq:force_balance_polynomial}
\vartwo{q,m}(\vartwo{F,m},z) + \vartwo{P,aux}\lambda_v \leq \vartwo{F,fc} + \vartwo{F,batt}.
\end{equation}

The convex polynomial $\vartwo{q,m}(\vartwo{F,m},z)$ can be guaranteed to accurately approximate $\vartwo{F,m}/\vartwo{\eta,m}(\vartwo{F,m},z)$ for all motors known, as efficiency is practically concave in power \citep{RN859} the reciprocal of which is convex \citep{boyd2004convex}. 

The remaining aspect to be covered is motor operational constraints. Motors operate within two regions depending on rotational velocity, a constant force region under the cutoff speed expressed by the simple bounds
\begin{equation}
\underline{\vartwo{F,m}} \leq\ \vartwo{F,m} \leq \overline{\vartwo{F,m}}
\end{equation}
and a constant power region above the cutoff speed expressed by the linear inequalities
\begin{subequations}
\begin{alignat}{1}
\underline{\vartwo{P,m}}\lambda_v & \leq \vartwo{F,m},\\ 
\vartwo{F,m} & \leq \overline{\vartwo{P,m}}\lambda_v.
\end{alignat}
\end{subequations}

\subsection{Fuel Cell}\label{sec:fuel_cell}

\begin{figure}[t!]
\centering
\includegraphics[width=8.4cm,keepaspectratio]{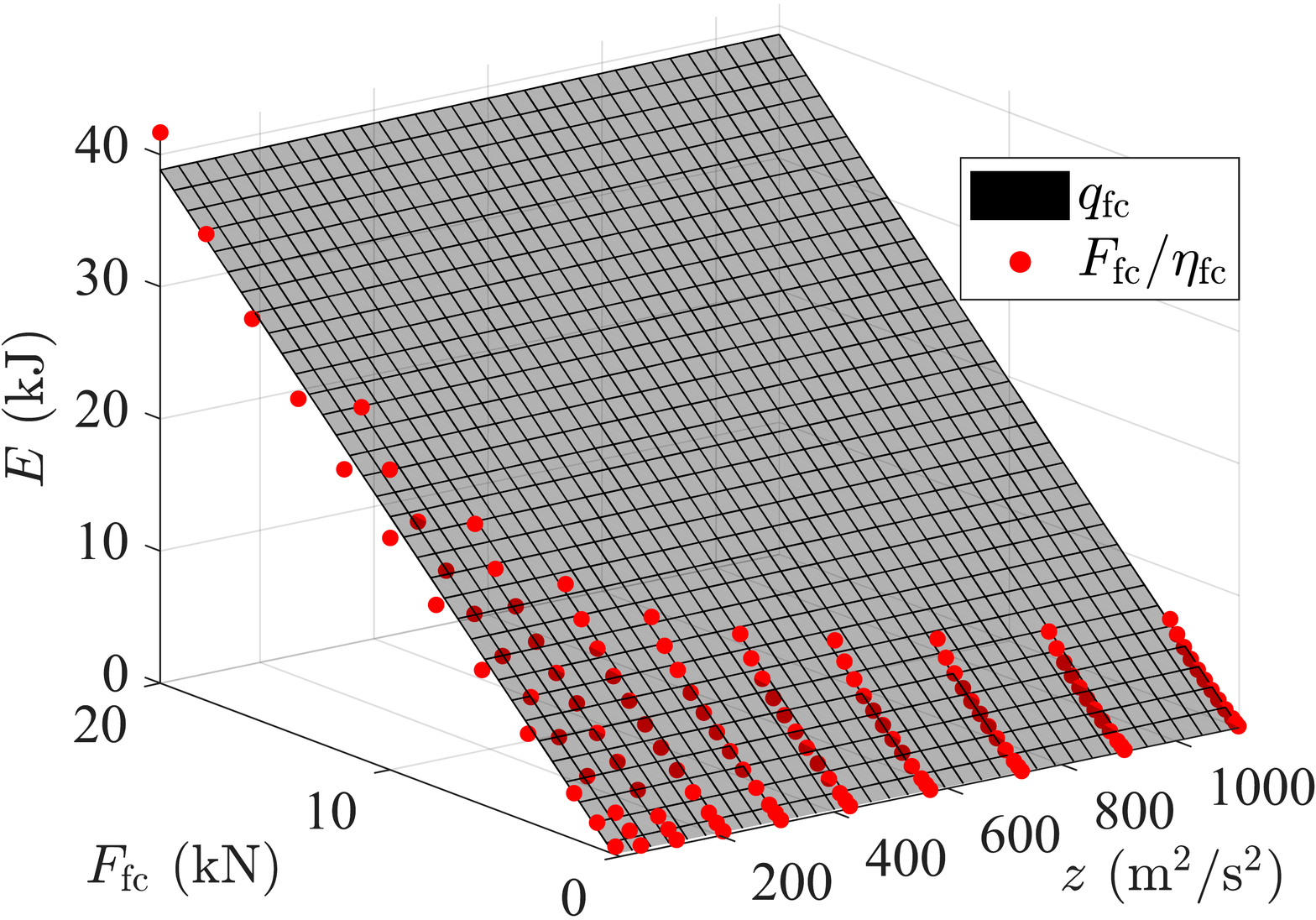}
\caption{Fuel cell data from \cite{RN879}. $\vartwo{\Delta,s}=1$.}
\label{fig:fc_consumption}
\end{figure}

To penalize hydrogen fuel consumption in the space-domain, we derive an expression for fuel energy consumed per longitudinal meter traveled. Using the look-up table efficiency model $\vartwo{\eta,fc}(\vartwo{F,fc},z)$, the exact fuel penalty per meter is $\vartwo{F,fc}/\vartwo{\eta,fc}(\vartwo{F,fc},z)$. \cite{RN779} proved using first-principle models that all fuel cell technologies admit a concave efficiency curve with power which implies that the exact penalty can be accurately approximated by the convex quadratic polynomial 
\begin{equation}\label{eq:q_fc}
\begin{split}
\vartwo{q,fc}(\vartwo{F,fc},z) := &p_{00}'+p_{10}'z+p_{01}'\vartwo{F,fc}\\
&+p_{11}'\vartwo{F,fc}v+p_{20}'z^2+p_{02}'\vartwo{F,fc}^2,
\end{split}
\end{equation}
as shown by Fig. \ref{fig:fc_consumption}.

The fuel cell power constraints are expressed by
\begin{subequations}
\begin{alignat}{1}
\underline{\vartwo{P,fc}}\lambda_v & \leq \vartwo{F,fc},\\ 
\vartwo{F,fc} & \leq \overline{\vartwo{P,fc}}\lambda_v,
\end{alignat}
\end{subequations}
where the lower bound $\underline{\vartwo{P,fc}}$ could be selected as strictly positive in order to curtail the excessive degradation that accompanies idling and restarting.

\subsection{Battery State-of-Charge}\label{sec:battery_soc}

Predicting the battery's state-of-charge $\zeta$ is vital in order to guarantee charge-sustaining operation---terminal battery charge identical to initial. The battery is modeled with a fixed open-circuit voltage $\vartwo{U,oc}$ and a fixed internal resistance $R$, a model that is accurate for the narrow state-of-charge range employed by hybrid vehicles. Experimental results by \cite{RN700} validate this model.

\cite{RN736} derive the change of state-of-charge
\begin{equation}\label{eq:dzeta_dt}
\Delta_\zeta = \frac{\vartwo{U,oc}-\sqrt{\vartwo{U,oc}^2-4\vartwo{P,batt}R}}{2R} \cdot \frac{1}{3600Q} \cdot \Delta_\text{t},
\end{equation}
where $Q$ is battery charge capacity, valid for $\vartwo{P,batt} \leq \vartwo{U,oc}^2/4R$. Accordingly, a positive/(negative) $\vartwo{P,batt}$ will discharge/(charge) the battery
\begin{equation}\label{eq:zeta}
\zeta_{i+1}=\zeta_i-\Delta_{\zeta,i}.
\end{equation}

\begin{figure}[t!]
\centering
\includegraphics[width=8.4cm,keepaspectratio]{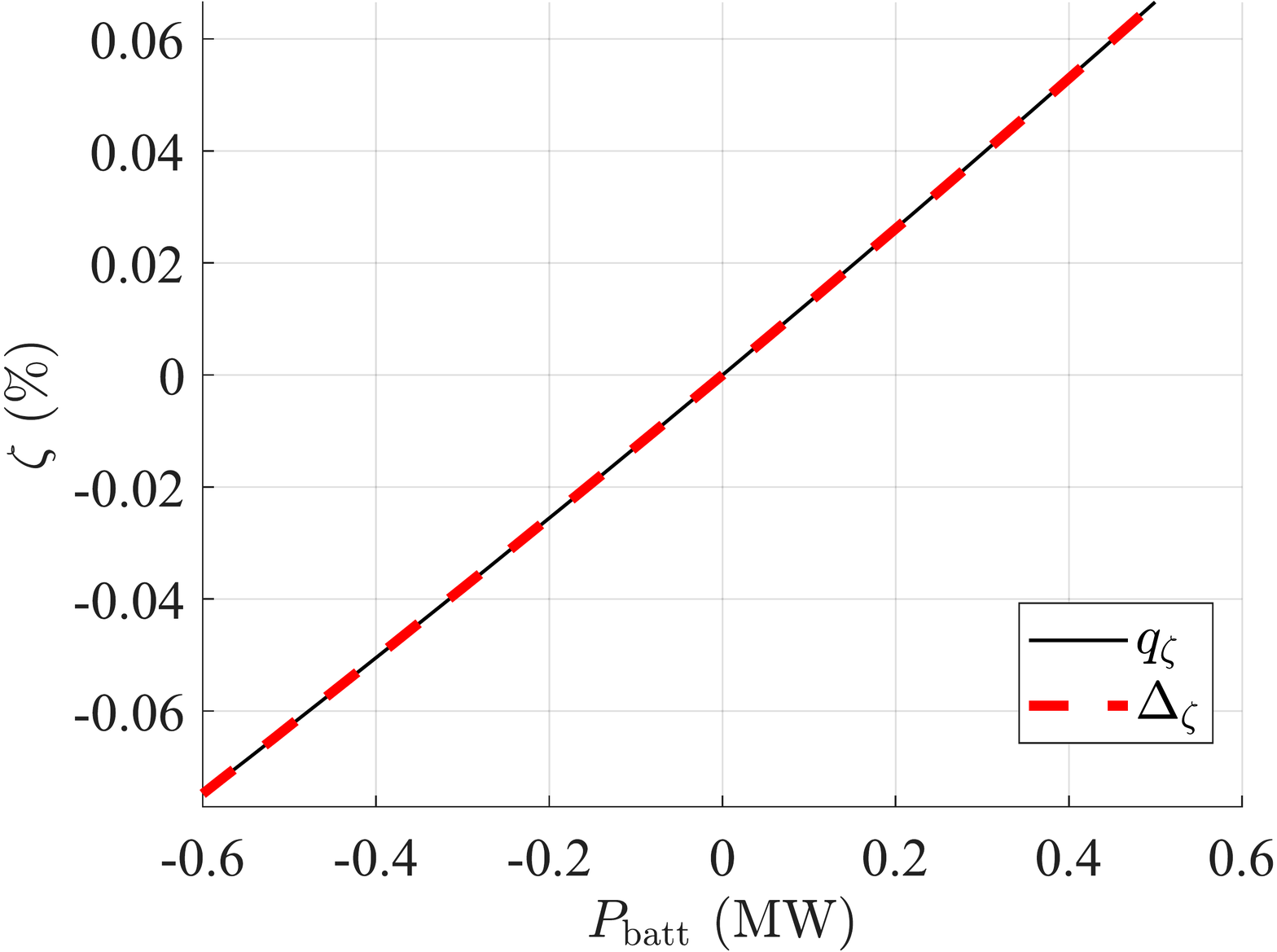}
\caption{Battery data from \cite{RN353}. $\vartwo{\Delta,t}=1$.}
\label{fig:batt_dzeta_dt}
\end{figure}

For a given $\vartwo{\Delta,t}$ the model \eqref{eq:dzeta_dt} is convex in $\vartwo{P,batt}$ because the negative sign preceding the concave square root flips it into a convex term. This empowers the convex quadratic polynomial $q_\zeta(\vartwo{P,batt}):=\alpha \vartwo{P,batt}^2 + \beta \vartwo{P,batt}$ to accurately approximate \eqref{eq:dzeta_dt}, as shown by Fig. \ref{fig:batt_dzeta_dt}. Nevertheless, an expression written in terms of spatial intervals $\vartwo{\Delta,s}$ rather than temporal $\vartwo{\Delta,t}$ needs to be found to complete a space-domain formulation. To derive such, start by assuming $\Delta_\zeta =q_\zeta(\vartwo{P,batt})\vartwo{\Delta,t}$
\begin{equation}
\Delta_\zeta = (\alpha\vartwo{P,batt}^2 + \beta \vartwo{P,batt})\vartwo{\Delta,t}
\end{equation}
which can be rewritten in terms of $\vartwo{F,batt}$ as
\begin{equation}\label{eq:F_batt_original}
\Delta_\zeta = (\alpha\vartwo{F,batt}^2v^2 + \beta \vartwo{F,batt}v)\vartwo{\Delta,t}
\end{equation}
followed by the substitution $v=\Delta_\text{s}/\Delta_\text{t}$
\begin{equation}
\Delta_\zeta = \Big(\alpha\vartwo{F,batt}^2v \frac{\Delta_\text{s}}{\Delta_\text{t}} + \beta \vartwo{F,batt} \frac{\Delta_\text{s}}{\Delta_\text{t}}\Big)\vartwo{\Delta,t}
\end{equation}
then cancel out $\Delta_\text{t}$ in order to obtain the spatial expression
\begin{equation}\label{eq:battery_model_1}
\Delta_\zeta = \alpha\vartwo{F,batt}^2v\vartwo{\Delta,s} + \beta \vartwo{F,batt}\vartwo{\Delta,s}.
\end{equation}

Equation \eqref{eq:battery_model_1} is non-convex but can be rewritten as
\begin{equation}\label{eq:battery_model_2}
\alpha\vartwo{F,batt}^2\vartwo{\Delta,s} = \frac{\Delta_\zeta-\beta \vartwo{F,batt}\vartwo{\Delta,s}}{v}
\end{equation}
then subsumed into
\begin{equation}\label{eq:battery_model_3}
\alpha\vartwo{F,batt}^2\vartwo{\Delta,s} = \lambda_\zeta\lambda_v
\end{equation}
using the linear auxiliary constraint
\begin{equation}\label{eq:lambda_zeta}
\lambda_\zeta=\Delta_\zeta-\beta \vartwo{F,batt}\vartwo{\Delta,s}
\end{equation}
and the convex constraint \eqref{eq:lambda_v}. 

The relaxation of the non-convex equality \eqref{eq:battery_model_3},
\begin{equation}\label{eq:battery_model_4}
\alpha\vartwo{F,batt}^2\vartwo{\Delta,s} \leq \lambda_\zeta\lambda_v,
\end{equation}
forms a convex feasible set for $\lambda_\zeta,\lambda_v \geq 0$ which is nonrestrictive, since $\lambda_v$ and the left-hand side of \eqref{eq:battery_model_4} are non-negative by definition.

\subsection{Battery Temperature}\label{sec:battery_temp}

Battery temperature $\vartwo{T,batt}$ is to be modeled in order to keep temperature under the upper bound
\begin{equation}\label{eq:t_batt_bound}
\vartwo{T,batt} \leq \overline{\vartwo{T,batt}}
\end{equation}
in order to preserve battery health. For a change in temperature of $\Delta\vartwo{T,batt}$ between intervals, battery temperature is predicted using the linear
\begin{equation}\label{eq:t_batt}
\varthree{T,batt,i+1}=\varthree{T,batt,i}+\Delta\varthree{T,batt,i}.
\end{equation}

Changes in temperature are caused by the electrochemical losses during use, the heat lost passively to the surroundings, and the heat extracted by the active cooling system. \cite{RN854} model the battery as a lumped mass $\vartwo{m,batt}$ with thermal capacity $\vartwo{c,batt}$ that admits a thermal content change of $\vartwo{m,batt}\vartwo{c,batt}\Delta\vartwo{T,batt}$ for a change $\Delta\vartwo{T,batt}$. Using the fictitious forces convention, the heat balance between spatial intervals is
\begin{equation}\label{eq:heat_balance_1}
\vartwo{m,batt}\vartwo{c,batt}\Delta\vartwo{T,batt} = (\vartwo{F,gen}-\vartwo{F,lost})\vartwo{\Delta,s},
\end{equation} 
where $\vartwo{F,gen}$ and $\vartwo{F,loss}$ denote the heat generated and lost per meter traveled, respectively.

\subsubsection{Derivation of Heat Generated}

$\vartwo{F,gen}$ can be directly expressed in terms of battery efficiency for both charging and discharging as
\begin{equation}\label{eq:f_gen_abs}
\vartwo{F,gen} = |\vartwo{F,batt}|\big(1-\vartwo{\eta,batt}(\vartwo{F,batt},v)\big),
\end{equation} 
where
\begin{equation}
\vartwo{\eta,batt}(\vartwo{F,batt},v):=
\begin{cases*}
\vartwo{U,batt}(\vartwo{F,batt}v)/\vartwo{U,oc}&$\vartwo{F,batt} \geq 0$ \\
\vartwo{U,oc}/\vartwo{U,batt}(\vartwo{F,batt}v)&otherwise
\end{cases*}
\end{equation}
and 
\begin{equation}
\vartwo{U,batt}(P):= \big(\vartwo{U,oc}+\sqrt{\vartwo{U,oc}^2-4\vartwo{P,batt}{R}}\big)/2.
\end{equation}

However, the equality \eqref{eq:heat_balance_1} cannot maintain its linear status if it were to admit the absolute value operation $|\vartwo{F,batt}|$ as required by \eqref{eq:f_gen_abs}. Alternatively, we propose to mimic $|\vartwo{F,batt}|$ by $\vartwo{F,dis}-\vartwo{F,chr}$ as in 
\begin{equation}\label{eq:f_gen_pos_neg}
\vartwo{F,gen} = (\vartwo{F,dis}-\vartwo{F,chr})\big(1-\vartwo{\eta,batt}(\vartwo{F,batt},v)\big),
\end{equation} 
where $\vartwo{F,dis}\geq\vartwo{F,batt},0$ and $\vartwo{F,chr}\leq\vartwo{F,batt},0$. Section 4 explains how $\vartwo{F,dis}$ and $\vartwo{F,chr}$ adopt the positive discharging and negative charging values of $\vartwo{F,batt}$, respectively. Lastly, the variable efficiency term $\vartwo{\eta,batt}(\vartwo{F,batt},v)$ in \eqref{eq:f_gen_pos_neg} impedes a linear expression due to its multiplication by the variables $\vartwo{F,dis}$ and $\vartwo{F,chr}$. Instead, we propose to simplify \eqref{eq:f_gen_pos_neg} using constant efficiency terms
\begin{equation}\label{eq:f_gen_pos_neg_eff}
\begin{split}
\vartwo{F,gen} = &\vartwo{F,dis}(1-\vartwo{\eta,dis})-\vartwo{F,chr}(1-\vartwo{\eta,chr}),
\end{split} 
\end{equation}
where $\vartwo{\eta,dis}$ and $\vartwo{\eta,chr}$ denote average discharging and charging efficiency, respectively.

\subsubsection{Derivation of Heat Lost}

The heat lost per meter traveled
\begin{equation}
\vartwo{F,lost}=\vartwo{F,amb}+\vartwo{F,act}
\end{equation}
comprises losses to ambient $\vartwo{F,amb}$ and active cooling system $\vartwo{F,act}$. The heat lost to ambient is easiest expressed as $h(\vartwo{T,batt}-\vartwo{T,amb})\vartwo{\Delta,t}$, where $h$ is rate of heat transfer per second. Upon substituting $\vartwo{\Delta,t}=1/v$ into the aforementioned
\begin{equation}\label{eq:f_amb}
\vartwo{F,amb} = h(\vartwo{T,batt}-\vartwo{T,amb})\lambda_v.
\end{equation}

\subsubsection{Compilation of Thermal Model}

Upon substituting and expanding \eqref{eq:f_gen_pos_neg_eff} and \eqref{eq:f_amb} into \eqref{eq:heat_balance_1} we get
\begin{equation}
\begin{split}
\vartwo{m,batt}\vartwo{c,batt}\Delta\vartwo{T,batt}= &\Big( \vartwo{F,dis}(1-\vartwo{\eta,dis})-\vartwo{F,chr}(1-\vartwo{\eta,chr})\\
&-h\vartwo{T,batt}\lambda_v+h\vartwo{T,amb}\lambda_v\\
&-\vartwo{F,act} \Big) \vartwo{\Delta,s}
\end{split}
\end{equation}
which is almost linear except for the term $h\vartwo{T,batt}\lambda_v$. Replace this non-linear term by the relaxed inequality
\begin{equation}\label{eq:lambda_T}
\lambda_T \leq h\vartwo{T,batt}\lambda_v
\end{equation}
for $\vartwo{T,batt},\lambda_v \geq 0$ to get the entirely linear
\begin{equation}\label{eq:compiled_thermal}
\begin{split}
\vartwo{m,batt}\vartwo{c,batt}\Delta\vartwo{T,batt}= &\Big( \vartwo{F,dis}(1-\vartwo{\eta,dis})-\vartwo{F,chr}(1-\vartwo{\eta,chr})\\
&-\lambda_T+h\vartwo{T,amb}\lambda_v\\
&-\vartwo{F,act} \Big) \vartwo{\Delta,s}.
\end{split}
\end{equation}

The non-negative condition imposed on $\vartwo{T,batt}$ is nonrestrictive, since a temperature of negative kelvin is physically infeasible. Section 4 explains how the inequality \eqref{eq:lambda_T} holds with equality when the upper temperature bound \eqref{eq:t_batt_bound} is approached.

\section{Optimization Formulation}

The models derived in sections 2 and 3 are now used to formulate the target optimization problem. The optimized system states are $(z,\zeta,\vartwo{T,batt})$; the main control variables are $(\vartwo{F,m},\vartwo{F,brk},\vartwo{F,fc},\vartwo{F,batt},\vartwo{F,act})$; and the auxiliary variables are $(v,\lambda_v,\lambda_\zeta,\lambda_T,\Delta_\zeta,\Delta\vartwo{T,batt},\vartwo{F,pos},\vartwo{F,neg})$. After obtaining the optimal trajectories to the fictitious force variables $(\vartwo{F,fc},\vartwo{F,batt},\vartwo{F,act})$, they are to be multiplied by the velocity trajectory in order to obtain their respective commands in terms of power.

The optimization problem computes the trajectory for $N$ intervals from $i=0,1,\cdots,N-1$ starting with initial states $(z_0,\zeta_0,\varthree{T,batt,0})$. The cost function
\begin{equation}\label{eq:cost_function}
\sum_i \Big(\vartwo{q,fc}(\varthree{F,fc,i},\vartwo{z,i}) + \varthree{F,act,i}\Big)\varthree{\Delta,s,i}
\end{equation}
penalizes hydrogen fuel consumption and the parasitic draw of the active cooling system.

The linear equality constraints \eqref{eq:kinetic_energy_2}, \eqref{eq:zeta}, and \eqref{eq:t_batt}, predict the system's states $(z,\zeta,\vartwo{T,batt})$, respectively. A second set of necessary linear equality constraints are \eqref{eq:lambda_zeta} and \eqref{eq:compiled_thermal} for the auxiliary variables $\lambda_\zeta$ and $\Delta\vartwo{T,batt}$. Moreover, the equality
\begin{equation}
\zeta_N = \zeta_0
\end{equation}
enforces charge-sustaining operation on the battery,
\begin{equation}\label{eq:journey_time}
\sum_i \varthree{\Delta,s,i}\lambda_{v,i} = \tau
\end{equation}
terminates the journey exactly $\tau$ seconds after start, and
\begin{equation}\label{eq:station_stops}
z_j = \vartwo{z,stop}
\end{equation}
halts the train at station stop intervals denoted $j$.

The linear inequality constraints are broken down into the simple lower and upper bounds
\begin{subequations}
\begin{alignat}{2}
\varthree{F,chr,i} & \leq \hspace{1em} 0 && \leq \lambda_{v,i},\lambda_{\zeta,i}, \lambda_{T,i}, \varthree{F,dis,i}\\
\underline{v} & \leq \hspace{0.9em} v_i && \leq \overline{v},\\ 
\underline{v}^2 & \leq \hspace{0.9em} z_i && \leq \overline{v}^2,\\ 
\underline{\zeta} & \leq \hspace{0.9em} \zeta_i && \leq \overline{\zeta},\\
0 & \leq \hspace{0em} \varthree{T,batt,i} && \leq \overline{\vartwo{T,batt}},\label{eq:temp_bound_inequality}\\
\underline{\vartwo{F,m}} & \leq \hspace{0.2em} \varthree{F,m,i} && \leq \overline{\vartwo{F,m}},\\ 
\underline{\vartwo{F,brk}} & \leq \hspace{0.1em} \varthree{F,brk,i} && \leq \overline{\vartwo{F,brk}}
\end{alignat}
\end{subequations}
and the more elaborate
\begin{subequations}
\begin{alignat}{2}
\underline{\vartwo{P,m}}\lambda_{v,i} & \leq \varthree{F,m,i} && \leq \overline{\vartwo{P,m}}\lambda_{v,i},\\ 
\underline{\vartwo{P,batt}}\lambda_{v,i} & \leq \varthree{F,batt,i} && \leq \overline{\vartwo{P,batt}}\lambda_{v,i},\\ 
\underline{\vartwo{P,fc}}\lambda_{v,i} & \leq \varthree{F,fc,i} && \leq \overline{\vartwo{P,fc}}\lambda_{v,i}.
\end{alignat}
\end{subequations}

Lastly, are the list of relaxed convex inequalities
\begin{subequations}\label{eq:relaxed_constraints}
\begin{alignat}{2}
1 & \leq v_i \lambda_{v,i},\label{eq:relaxed_1}\\
v_i^2 & \leq z_i,\label{eq:relaxed_2}\\
\vartwo{q,m}(\varthree{F,m,i},z_i) + \varthree{P,aux,i}\lambda_{v,i} & \leq \varthree{F,fc,i} + \varthree{F,batt,i},\label{eq:relaxed_3}\\
\alpha \varthree{F,batt,i}^2\varthree{\Delta,s,i} & \leq \lambda_{\zeta,i} \lambda_{v,i},\label{eq:relaxed_4}\\
\lambda_{T,i} & \leq h\varthree{T,batt,i}\lambda_{v,i},\label{eq:relaxed_5}\\
\varthree{F,chr,i} & \leq \varthree{F,batt,i},\label{eq:relaxed_6}\\
\varthree{F,batt,i} & \leq \varthree{F,dis,i}.\label{eq:relaxed_7}
\end{alignat}
\end{subequations}

The constraint \eqref{eq:relaxed_1} implies that $v$ is strictly positive and thus $z$ as well due to \eqref{eq:relaxed_2}. Nevertheless, in order to emulate being stationary at station stops in \eqref{eq:station_stops}, $\vartwo{z,stop}$ is set to a small positive value that approaches zero. During station stops $\varthree{F,ext,j}$ is zeroed in order to successfully emulate a stationary state with brakes locked (see \eqref{eq:kinetic_energy_2}). Since the optimized speed profile is strictly positive, the sampling intervals during station stops $\Delta_{\text{s},j}$ are adjusted \textit{a priori} to the multiplication of dwell (wait) time by $\sqrt{\vartwo{z,stop}}$. Although the optimized speed at station stops never attains zero, in practice, it can be zeroed without affecting feasibility or optimality if $\vartwo{z,stop}$ was small enough.

In order to prove the optimality of the proposed formulation, all relaxed constraints \eqref{eq:relaxed_constraints} need to be proven to hold with equality. The following justifies inequality tightness:

\begin{itemize}
\item \eqref{eq:relaxed_1}: the summation $\sum_i \lambda_{v,i}$ is fixed through \eqref{eq:journey_time} and $v$ has the incentive to drop due to losses in \eqref{eq:f_ext};
\item \eqref{eq:relaxed_2}: $z$ has incentive to drop due to penalty \eqref{eq:cost_function} but $v$ is constrained from beneath by \eqref{eq:relaxed_1};
\item \eqref{eq:relaxed_3}: $\vartwo{F,batt}$ has incentive to go negative to gather free charge and minimize $\vartwo{F,fc}$ while $\vartwo{q,m}$ has incentive to move the train to fulfill journey time \eqref{eq:journey_time};
\item \eqref{eq:relaxed_4}: the original expression \eqref{eq:F_batt_original} when relaxed, $\Delta_\zeta \geq (\alpha\vartwo{F,batt}^2v^2 + \beta \vartwo{F,batt}v)\vartwo{\Delta,t}$, would rather have positive $\vartwo{F,batt}$ to move the train and push $\Delta_\zeta$ to zero or negative to gain free charge;
\item \eqref{eq:relaxed_5},\eqref{eq:relaxed_6},\eqref{eq:relaxed_7}: if the upper temperature bound in \eqref{eq:temp_bound_inequality} is reached, \eqref{eq:compiled_thermal} would rather tighten \eqref{eq:relaxed_5}, \eqref{eq:relaxed_6}, and \eqref{eq:relaxed_7}, before relying on the active cooling system command $\vartwo{F,act}$ that is penalized in \eqref{eq:cost_function}.
\end{itemize}

The optimization problem proposed above is convex because it penalizes a convex quadratic cost function subject to linear equality and convex inequality constraints. It can be formulated and solved as a second-order cone program.

\section{Conclusion}
Models for the the concurrent optimization of hybrid train speed, EMS, and battery thermals, were presented. A relaxed convex problem was formulated in order to alleviate computational concerns while the tightness of the relaxed constraints was justified. The accuracy of the proposed convex models was proven by graphical means and analyzing the convexity properties of the original first-principle models. The benefit from this holistic optimization approach is yet to be verified on a real case study, after which optimizing fuel cell thermals and optimizing the operation of singular fuel cell stacks independently is to be investigated.

\bibliography{ifacconf}             

\end{document}